\documentclass[12pt]{article}
 \usepackage{rotating}
 \usepackage[export]{adjustbox}
 \usepackage{bm}
 \usepackage{empheq}
 \usepackage[most]{tcolorbox}
 \usepackage{caption}
 \usepackage{graphicx}
 \usepackage{subfig}
 \usepackage{wasysym}
 \usepackage[margin=1in]{geometry} 
 \usepackage{amsmath,amsthm,amssymb}
 \usepackage{hyperref}
 \usepackage{algorithm2e}
 \usepackage{algorithmic}
 \usepackage{float}
 \usepackage{multimedia}
 \usepackage{tikz} 
 \usepackage{atbegshi,picture}

 \newenvironment{theorem}[2][Theorem]{\begin{trivlist}
 		\item[\hskip \labelsep {\bfseries #1}\hskip \labelsep {\bfseries 
 		#2.}]}{\end{trivlist}}

\title{A short proof that sweeping is always 
possible for a spatial discretization with regular triangles and no hanging 
nodes}
\author{Thomas Camminady\thanks{Corresponding author: 
\texttt{camminady@kit.edu}}, 
Martin Frank\\
	Karlsruhe Institute of Technology}
\begin{document}

\maketitle

\abstract{Sweeping is a commonly used procedure to explicitly solve the 
discrete ordinates equation, which itself is a common approximation of the 
neutron transport equation. To sweep through the computational domain, an 
ordering of the spatial cells is required that obeys the flow of information. 
We show that this ordering can always be found, assuming a discretization of 
the spatial domain with regular triangles with no hanging nodes. 
}

\section{Introduction}
\begin{figure} 
	\centering
\begin{minipage}{0.3\textwidth}
	\includegraphics[width=0.9\linewidth]{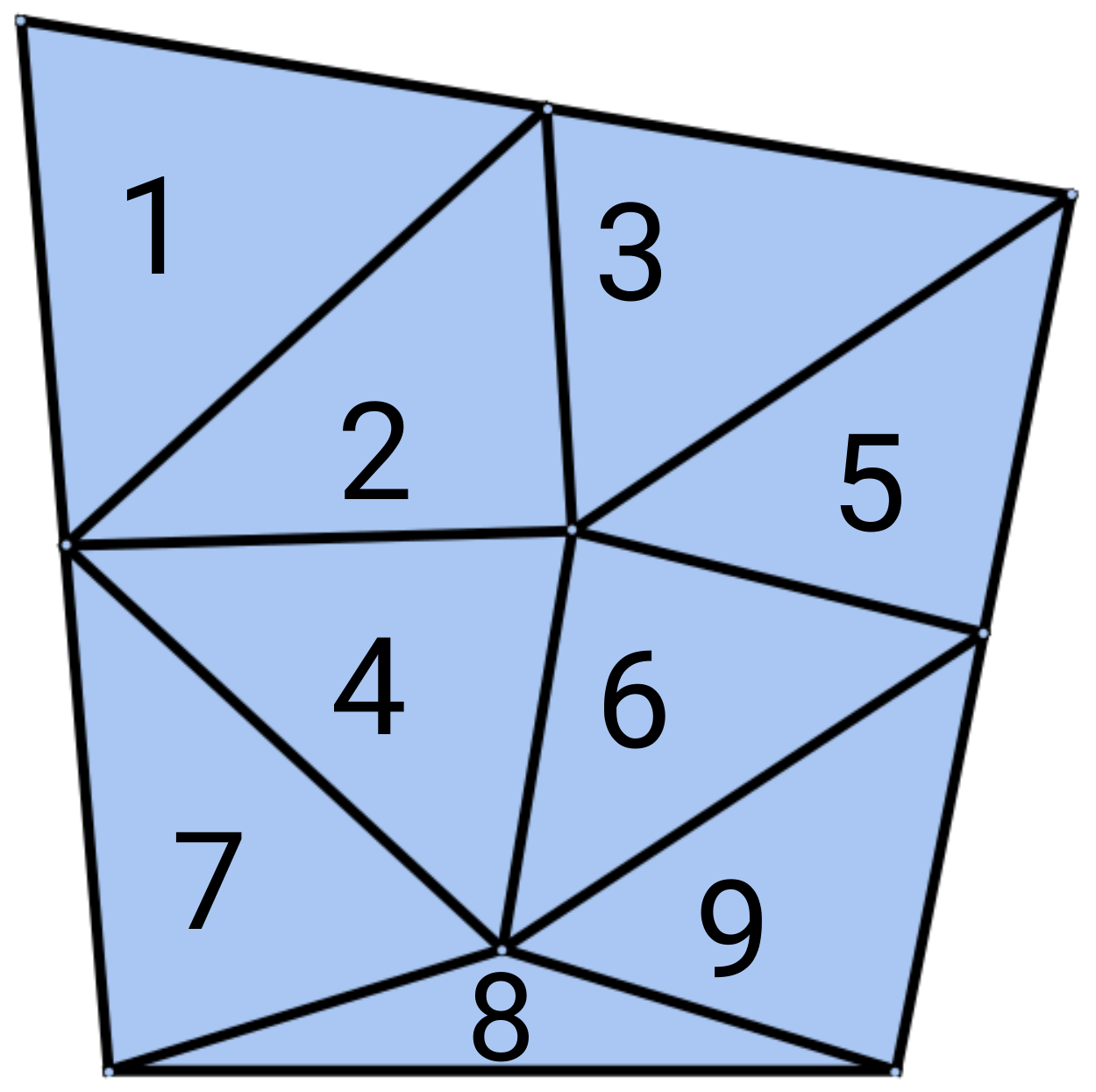}
	\caption{Triangulation of the computational domain with labels.}
	\label{fig:smallmeshlabels}
\end{minipage}
\hspace{0.5cm}
\begin{minipage}{0.3\textwidth}
	\includegraphics[width=0.9\linewidth]{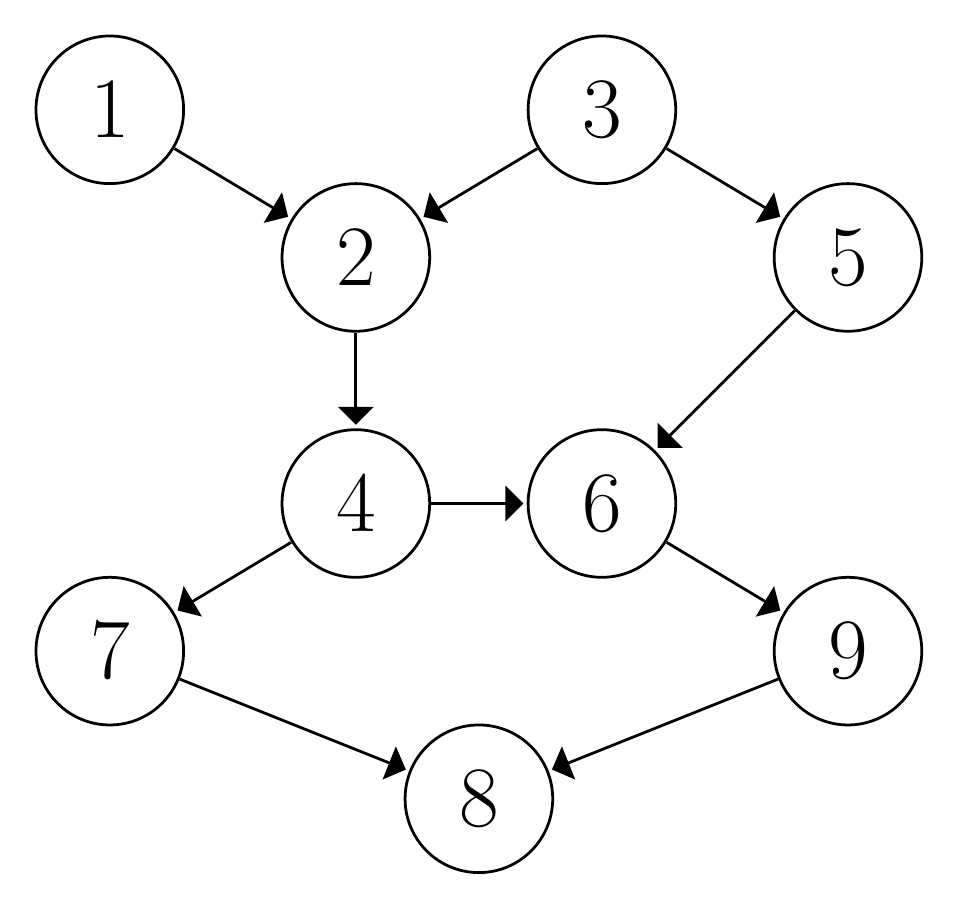}
	\caption{Induced dependency graph for \\$\Omega=(0,-1)^T$.}
	\label{fig:directegraph1}
\end{minipage}
\hspace{0.5cm}
\begin{minipage}{0.3\textwidth}
	\includegraphics[width=0.9\linewidth]{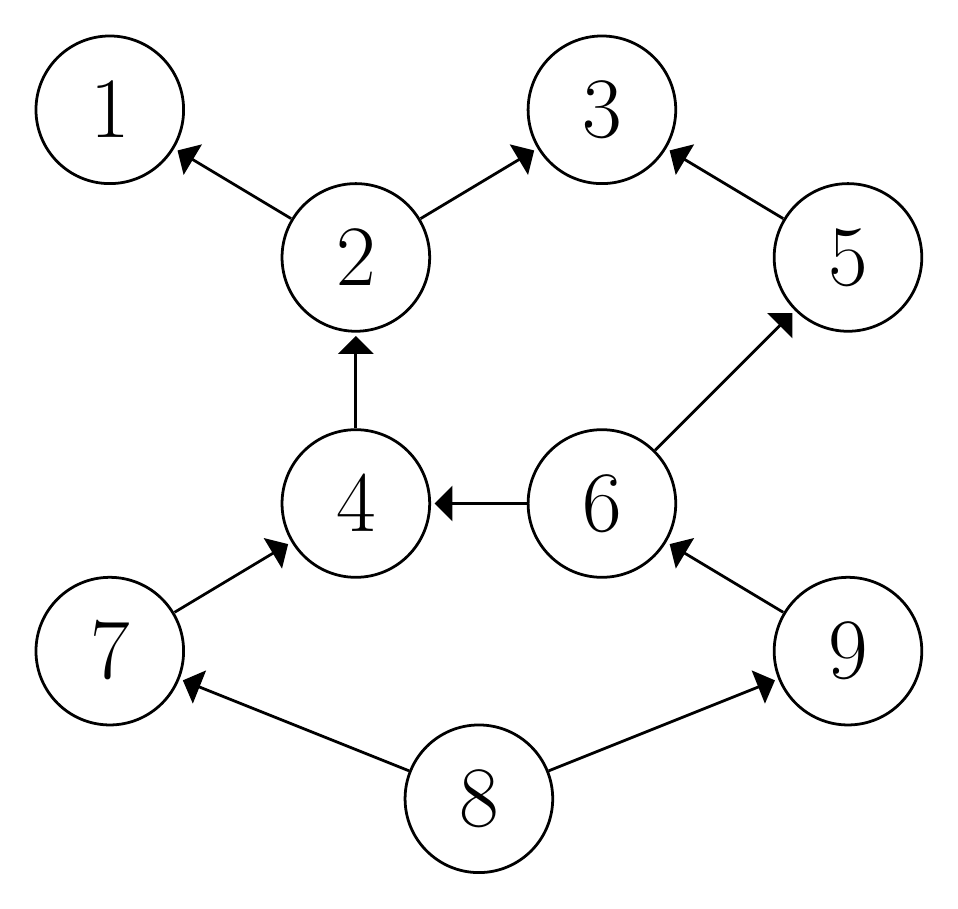}
	\caption{Induced dependency graph for \\$\Omega=(0,+1)^T$.}
	\label{fig:directegraph2}
\end{minipage} 
\end{figure}

\begin{figure}[h!]%
	\centering
	\includegraphics[width=0.3\linewidth]{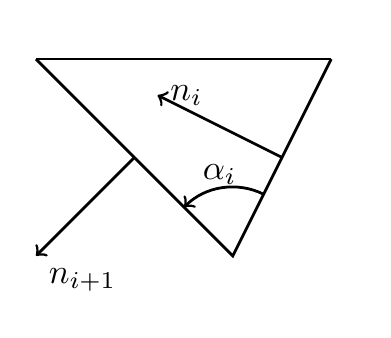}
	\caption{Incoming and outgoing normal vector for cell $C_i$ and the 
		corresponding angle $\alpha_i$.}
	\label{fig:cellci}
\end{figure}

 Iterative methods for the numerical solution of transport processes often 
 make use of sweeping to invert the streaming operator \cite{adams2002fast}.
 Sweeping requires to march through all spatial cells in a way that obeys the 
 flow of information, prescribed by a direction $\Omega$.
 By that, a dependency between cells is induced which implies a dependency 
 graph $G$.
 Here, two nodes $i'$ and $i$ share a directed edge, pointing from 
 $i'$ to $i$, if spatial cell $C_i$ does depend on spatial cell $C_{i'}$.
 To be more precise: $C_i$ depends on $C_{i'}$, if and only if they share 
 an edge $e_{i',i}$ and the 
 inward pointing normal $n_i$ (inward pointing with respect to cell $C_i$) 
 onto edge $e_{i',i}$ satisfies $\langle n_i,\Omega\rangle >0$. Fig. 
 \ref{fig:smallmeshlabels} shows a triangulation with labeled cells. For two 
 different directions $\Omega$, Fig. \ref{fig:directegraph1} and 
 \ref{fig:directegraph2} show the induced dependency for $\Omega$ downward 
 pointing and $\Omega$ upward pointing, respectively.
The inward and 
outward pointing normal vectors are sketched for a cell $C_i$ in Fig. 
\ref{fig:cellci}.
 
 Every directed acyclic graph has a unique topological sorting 
 \cite{skiena2008section}. This
 allows to march through the cells in a proper way, i.e. visiting cell $C_i$ 
 after visiting all cells that it depends on.
 We will show that the 
 dependency graph $G$ is indeed acyclic for a triangulation of regular 
 triangles and no hanging nodes and sweeping 
 is therefore always possible. 

While a lot of research has focused onto the development of fast sweeping 
strategies, especially for parallel architectures 
\cite{bailey2008analysis,baker1998sn,koch1992parallel,koch1992solution}, the 
authors are not aware of a rigorous proof that sweeping is always possible 
under 
certain 
conditions.
 
% Fig. \ref{fig:discretizationAndOrdering} shows two resulting dependency 
%graphs 
% for the same triangulation, each for a different direction $\Omega$.

\newpage

\section{Main Result}
\begin{figure}
\centering
\begin{minipage}{0.3\textwidth}
	\includegraphics[width=1\linewidth,left]{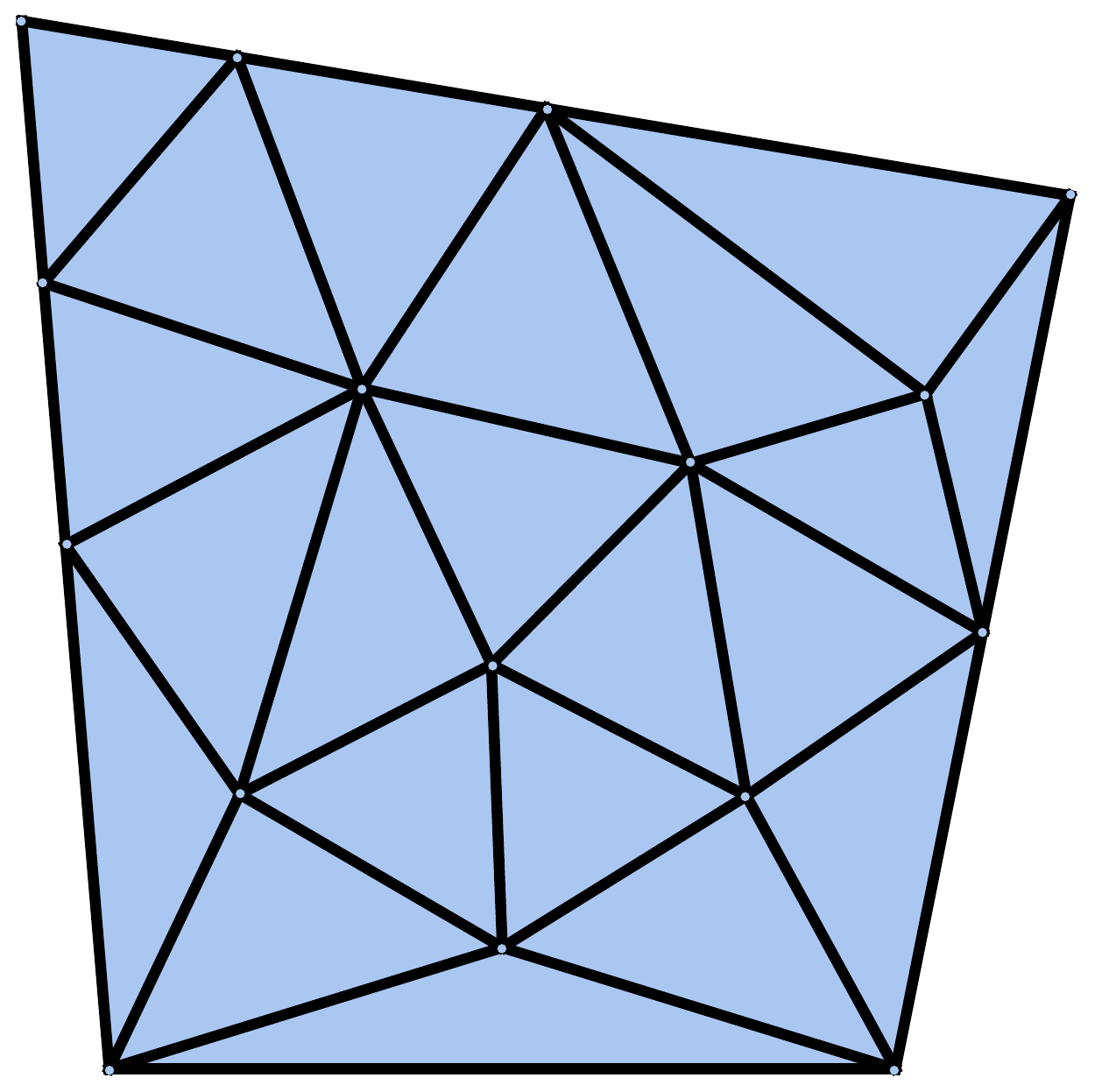}
	\caption{Triangulation of an arbitrary domain.}
	\label{fig:mesh}
\end{minipage}
\hspace{3cm}
\begin{minipage}{0.3\textwidth}
	\includegraphics[width=1\linewidth,right]{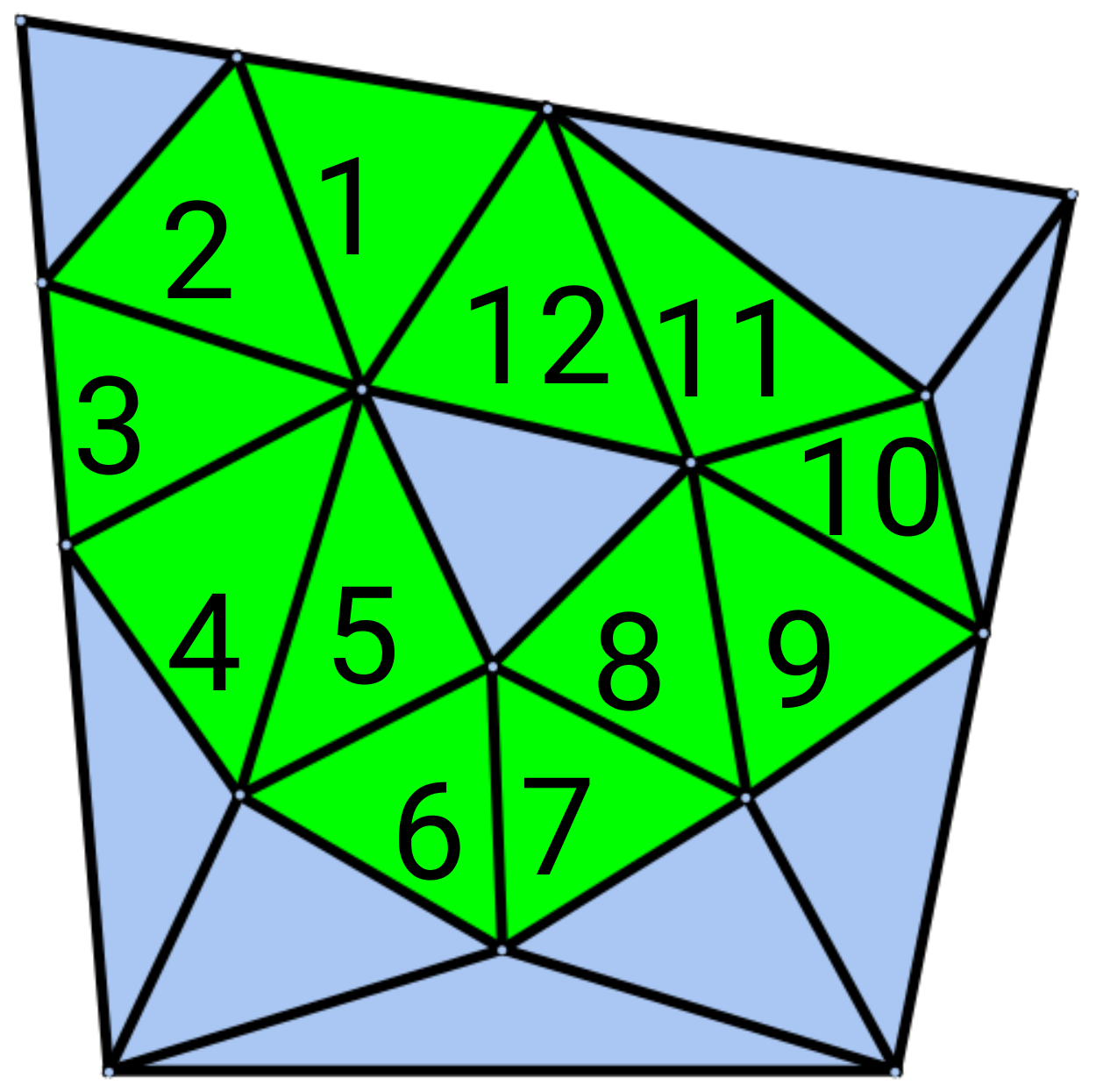}
	\caption{A cycle within the given triangulation.}
	\label{fig:cycle}
\end{minipage}
\end{figure}

\begin{theorem}{1}	
Consider a domain that is discretized by a set of cells 
$\{C_i\}_{i=1,\dots,I}$ where each cell $C_i$ is a triangle and we do not allow 
for hanging nodes, shown in Fig. \ref{fig:mesh}. 
Furthermore, we have a 
fixed direction $\Omega$ that 
prescribes the flow of 
information. Under these conditions, sweeping is possible.
\end{theorem}

\begin{figure}[]
\centering
\includegraphics[width=0.3\linewidth]{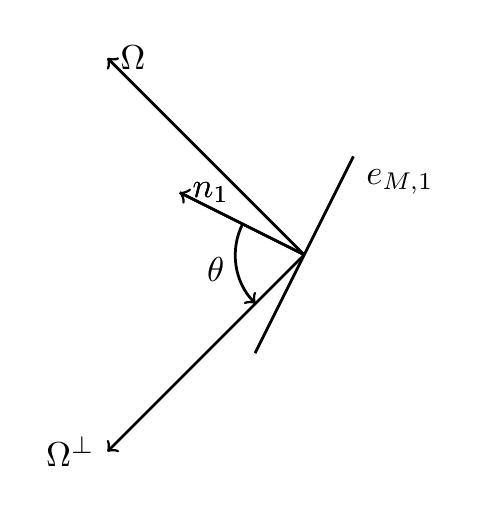}
\caption{The normal $n_1$, with $\Omega$ and $\Omega^\perp$, as well as the 
	angle $\theta$ between $n_1$ and $\Omega^\perp$.}
\label{fig:vectors2}
\end{figure}

\begin{proof}
We know that sweeping is possible, if and only if the induced dependency 
graph 
$G$ is acyclic. 
Now assume the dependency graph is not acyclic and provoke a 
contradiction.

If the dependency graph $G$ is not acyclic, there exists a sequence of 
cells \\
$\{C_{i_1},C_{i_2}, \dots,C_{i_M},C_{i_1}\}$, such that two consecutive 
cells 
share an edge 
$e_{{i_j},{i_{j+1}}}$ and $\langle n_{i_j},\Omega_k\rangle >0, \forall 
j=1,\dots,M$. Without loss of generality, let $i_j=j$ and assume that we 
pass 
through the cells in a counterclockwise manner as shown in Fig 
\ref{fig:cycle}.

Now label the angles between $n_i$ and $n_{i+1}$ by 
$\alpha_{i}$ as shown in Fig.\ref{fig:cellci}. A counterclockwise turn 
corresponds to $\alpha_{i}>0$ and 
a 
clockwise turn to $\alpha_{i}<0$. 
We know that $-\pi<\alpha_{i}<\pi$ as we 
consider regular triangles. Since we perform a full counterclockwise turn, 
$\sum_{i=1}^M \alpha_i=2\pi$.
Denote the angle between $n_1$ and $\Omega^\perp$ by $\theta$ as sketched 
in 
Fig. \ref{fig:vectors2}, with $0<\theta<\pi$ and 
$\Omega^\perp$ the vector normal to 
$\Omega$.
Let $R_{\alpha_i}$ be the rotation matrix that encodes rotating with 
magnitude $\alpha_i$ around the $z$-axis.
Then 
\begin{align*}
	n_{i+1} =R_{\alpha_{i}}n_i = \prod_{k=1}^{i} R_{\alpha_{k}} 
	n_1= R_{\sum_{k=1}^{i} \alpha_k} n_1.
\end{align*}	
If we turn $n_1$ (counterclockwise) by more than $\theta$ but less than 
$\theta+\pi$, then $\langle n_1,\Omega\rangle<0$.
However, there exist $i^*$ such that 
$\sum_{k=1}^{i^*-1} \alpha_k \leq \theta$, but $\theta<\sum_{k=1}^{{i^*}} 
\alpha_k < 
\theta+\pi<2\pi$, as $-\pi<\alpha_{i^*}<\pi$.
Then $n_{i^*+1}$ is $n_1$ turned (counterclockwise) by more than $\theta$, 
but 
less 
than $\theta+\pi$.
Therefore $\langle n_{i^*+1}, \Omega\rangle <0$ which contradicts the 
assumption 
and finishes 
the proof.
\end{proof}

\begin{figure}
	\centering
	\includegraphics[width=0.7\linewidth]{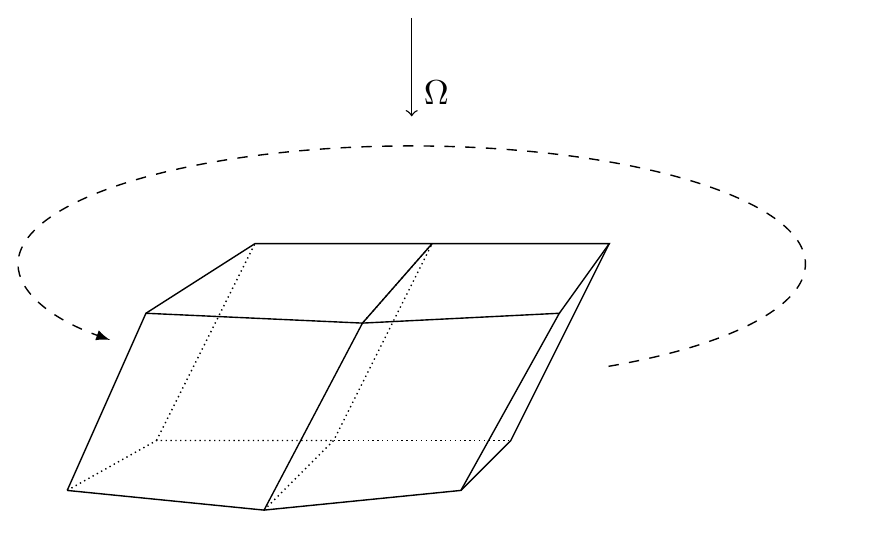}
	\caption{Quadrilaterals can be arranged in a circle such that no 
		topological 
		ordering can be obtained for $\Omega=(0,0,-1)^T$. }
	\label{fig:parallelogram}
\end{figure}

\section{The three dimensional case}
The analogue setup for the three dimensional case would be a triangulation with 
tetrahedra, again with no hanging nodes.
% Here, the previously described chain 
%of arguments does no longer hold, since we have rotations around different 
%axes 
%(always the edge that three consecutive tetrahedra share). Rewriting the 
%product 
%of rotation matrices as the rotation matrix of the sum of the angles is 
%therefore 
%no longer possible.
In general, allowing arbitrary convex quadrilaterals for the triangulation of a 
computational domain does not imply the absence of circular dependencies. This 
is sketched in Fig. \ref{fig:parallelogram}. Quadrilaterals can be arranged in 
a plane, perpendicular to the $z$-axis such that they form a circular 
dependency for a direction of flow along the $z$-axis.

For tetrahedra however, the authors assume the analogue result to hold as for 
the two dimensional case. In the absence of hanging nodes, sweeping should be 
possible. Turning the previous proof into a three dimensional version appears 
difficult, as rotation matrices no longer rotate around the same axis. 
Therefore different techniques have to be investigated to proof this claim.
\bibliographystyle{plain}
\bibliography{Camminady}

\begin{thebibliography}{1}

\bibitem{adams2002fast}
Marvin~L Adams and Edward~W Larsen.
\newblock Fast iterative methods for discrete-ordinates particle transport
  calculations.
\newblock {\em Progress in nuclear energy}, 40(1):3--159, 2002.

\bibitem{bailey2008analysis}
Teresa~S Bailey and Robert~D Falgout.
\newblock Analysis of massively parallel discrete-ordinates transport sweep
  algorithms with collisions.
\newblock Technical report, Lawrence Livermore National Laboratory (LLNL),
  Livermore, CA, 2008.

\bibitem{baker1998sn}
Randal~S Baker and Kenneth~R Koch.
\newblock An sn algorithm for the massively parallel cm-200 computer.
\newblock {\em Nuclear Science and Engineering}, 128(3):312--320, 1998.

\bibitem{koch1992parallel}
Kenneth~R Koch, Randal~S Baker, and Raymond~E Alcouffe.
\newblock A parallel algorithm for 3d sn transport sweeps.
\newblock {\em Technical Report}, 1992.

\bibitem{koch1992solution}
Kenneth~R Koch, Randal~S Baker, and Raymond~E Alcouffe.
\newblock Solution of the first-order form of the 3-d discrete ordinates
  equation on a massively parallel processor.
\newblock {\em Transactions of the American Nuclear Society}, 65(108):198--199,
  1992.

\bibitem{skiena2008section}
Steven~S Skiena.
\newblock Section 15.2: Topological sorting.
\newblock {\em The Algorithm Design Manual (2nd ed.), Springer-Verlag, London},
  pages 481--483, 2008.

\end{thebibliography}

\end{document}